\documentclass[11pt,a4]{elsarticle}

\usepackage{setspace}
\usepackage{geometry}
\usepackage{url}
\usepackage{graphicx,graphics}
\usepackage{subfigure}
\usepackage{epsfig}
\usepackage{psfrag}
\usepackage{amsmath,amssymb,amsfonts,amsthm}
\usepackage{mathrsfs}
\usepackage{color}
\usepackage{float}
\usepackage{natbib}
\usepackage{tabularx}
\usepackage{multirow}

\newcommand{\Eref}[1]{Equation (\ref{#1})}
\newcommand{\fref}[1]{Figure (\ref{#1})}

\newcommand{\bveps}{\boldsymbol{\varepsilon}}
\newcommand{\bvsig}{\boldsymbol{\sigma}}

\newcommand{\cc}{\mathbf{C}}
\newcommand{\dd}{\mathbf{D}}

\newcommand{\kk}{\mathbf{K}}

\newcommand{\qq}{\mathbf{q}}

\newcommand{\uu}{\mathbf{u}}

\begin{document}

\begin{frontmatter}


\title{Static bending and free vibration of cross ply laminated composite plates using NURBS based finite element method and unified formulation}

\author[unsw]{S Natarajan \corref{cor1}\fnref{fn1}}
\author[viet]{Hung Nguyen-Xuan}
\author[portugal,dubai]{AJM Ferreira}
\author[italy,dubai]{E Carrera}


\address[unsw]{School of Civil \& Environmental Engineering, The University of New South Wales, Sydney, Australia}
\address[viet]{Department of Mechanics, Faculty of Mathematics \& Computer Science, University of Science, HCMC, Vietnam}
\address[portugal]{Faculdade de Engenharia da Universidade do Porto, Porto, Portugal.}
\address[italy]{Department of Aeronautics and Aerospace Engineering, Politecnico di Torino, Italy.}
\address[dubai]{Department of Mathematics, Faculty of Science, King Abdulaziz University, P.O. Box 80203, Jeddah 21589, Saudi Arabia.}

\begin{abstract}
This paper presents an effective formulation to study the response of laminated composites based on isogeometric approach (IGA) and Carrera unified formulation (CUF). The IGA utilizes the non-uniform rational B-spline (NURBS) functions which allows to construct higher order smooth functions with less computational effort. The static bending and the free vibration of thin and moderately thick laminates plates are studied. The present approach also suffers from shear locking when lower order functions are employed and the shear locking is suppressed by introducing a modification factor. The combination of the IGA with the CUF allows a very accurate prediction of the field variables. The effectiveness of the formulation is demonstrated through numerical examples.
\end{abstract}


\end{frontmatter}


\section{Introduction} The need for high strength-to and high stiffness-to-weight ratio materials has led to the development of laminated composite materials. This class of material has seen increasing utilization as structural elements, because of the possibility to tailor the properties to optimize the structural response. Since, its inception, different approaches have been employed to study the response of such laminated composite plates, ranging from complete 3D analysis to two-dimensional theories. A brief overview on the development of different plate theories is given in~\cite{khandannoroozi2012,mallikarjunakant1993}. The various two dimensional plate theories can be further classified into three different approaches: (a) equivalent single layer theories~\cite{reddy1984}; (b) discrete layer theories~\cite{guonagy2014} and (c) mixed plate theory. Among these, the equivalent single layer theories, viz., first order shear deformation theory~\cite{rolfesrohwer1997}, second and higher order accurate theory~\cite{reddy1984,kantswaminathan2001} are the most popular theories employed to describe the plate kinematics. Existing approaches in the literature to study plate and shell structures made up of laminated composites uses the finite element method based on Lagrange basis functions~\cite{ganapathipolit1996} or non-uniform rational B splines (NURBS)~\cite{kapoorkapania2012} or meshfree methods~\cite{liewzhao2011}. Not only these approaches suffer from shear locking when applied to thin plates, these techniques does not provide a single platform to test the performance of various theories. Thanks to the recent derivation of series of axiomatic approaches by Carrera~\cite{carrera2001}, coined as Carrera Unified Formulation~\cite{carrerademasi2002} for the general description of two-dimensional formulations for multilayered plates and shells. With this unified formulation, it is possible to implement in a single software a series of hierarchical formulations, thus affording a systematic assessment of different theories ranging from simple equivalent single layer models up to higher order layerwise descriptions.

Recent interest in the unified formulation has led to the development of discrete models such as those based on finite element method~\cite{carreracinefra2010,natarajanferreira2013} and more recently meshless methods~\cite{liewzhao2011}. Nevertheless, even with the unified framework, there is an important shortcoming. With lower order basis functions within the finite element framework, when applied to thin plates, the formulation suffers from shear locking. Intensive research over the the past decades has to led to some of the robust methods to suppress the shear locking syndrome. This includes: (a) reduced integration~\cite{hughescohen1978}; (b) use of assumed strain method~\cite{nguyen-xuanrabczuk2008}; (c) using field redistributed shape functions~\cite{somashekarprathap1987}; (d) mixed interpolation tensorial components (MITC) technique with strain smoothing~\cite{bathedvorkin1985} and (e) very recently, twist Kirchhoff plate element~\cite{santosevans2012} and the 3D consistent formulation based on the scaled boundary finite element method~\cite{mansong2013}.

The main objective of this manuscript is to investigate the potential application of the NURBS based isogeometric finite element method within the Carrera Unified Formulation (CUF) to study the global response of cross-ply laminated composites. The present formulation also suffers from shear locking when lower order basis functions are employed to thin plates. To address lower-order NURBS elements for plates, the introduction of a stabilization technique into shear locking has been studied in \cite{ChienHungNhon2012}. The other approach to suppress is to employ higher order basis functions~\cite{veigabuffa2012}. In this study, to alleviate shear locking, a simple modification is done to the shear term when lower order NURBS basis functions are used. However, the draw back of this approach is that the shear correction factor is problem dependent. The influence of various parameters, viz., the ply thickness, the ply orientation, the plate geometry, the material property and the boundary conditions on the global response is numerically studied.

The paper commences with a brief discussion on the unified formulation for plates and the finite element discretization. Section \ref{isog} describes the isogeometric approach employed in this study, followed by a technique to address shear locking when lower order NURBS functions are used to discretize the field variables. The efficiency of the present formulation, numerical results and parametric studies are presented in Section~\ref{numres}, followed by concluding remarks in the last section.

\section{Carrera Unified Formulation}\label{cuftheory}
\subsection{Basis of CUF}
Let us consider a laminated plate composed of perfectly bonded layers with coordinates $x,y$ along the in-plane directions and $z$ along the thickness direction of the whole plate, while $z_k$ is the thickness of the $k^{\rm th}$ layer. The CUF is a useful tool to implement a large number of two-dimensional models with the description at the layer level as the starting point. By following the axiomatic modelling approach, the displacements $\uu(x,y,z) = ( u(x,y,z), v(x,y,z), w(x,y,z))$ are written according to the general expansion as:
\begin{equation}
\uu(x,y,z) = \sum\limits_{\tau = 0}^N F_\tau(z) \uu_\tau(x,y)
\label{eqn:unifieddisp}
\end{equation}
where $F(z)$ are known functions to model the thickness distribution of the unknowns, $N$ is the order of the expansion assumed for the through-thickness behaviour. By varying the free parameter $N$, a \emph{hierarchical} series of two-dimensional models can be obtained. The strains are related to the displacement field via the geometrical relations:
\begin{eqnarray}
\bveps_{pG} = \left[ \begin{array}{ccc} \varepsilon_{xx} & \varepsilon_{yy} & \gamma_{xy}  \end{array} \right]^{\rm T} = \dd_p \uu \nonumber \\
\bveps_{nG} = \left[ \begin{array}{ccc} \gamma_{xz} & \gamma_{yz} & \varepsilon_{zz} \end{array} \right]^{\rm T} = \left( \dd_{np} + \dd_{nz} \right) \uu
\label{eqn:strainDisp}
\end{eqnarray}
where the subscript $G$ indicate the geometrical equations, $\dd_{p}, \dd_{np}$ and $\dd_{nz}$ are differential operators given by:
\begin{eqnarray}
\dd_p = \left[ \begin{array}{ccc} \partial_x & 0 & 0 \\ 0 & \partial_y & 0 \\ \partial_y & \partial_x & 0 \end{array} \right], \hspace{0.5cm} \dd_{np} = \left[ \begin{array}{ccc} 0 & 0 & \partial_x \\ 0 & 0 & \partial_y \\ 0 & 0 & 0 \end{array} \right], \nonumber \\
\dd_{nz} = \left[ \begin{array}{ccc} \partial_z & 0 & 0 \\ 0 & \partial_z & 0 \\ 0 & 0 & \partial_z \end{array} \right].
\end{eqnarray}
The 3D constitutive equations are given as:
\begin{eqnarray}
\bvsig_{pC} = \cc_{pp} \bveps_{pG} + \cc_{pn} \bveps_{nG} \nonumber \\
\bvsig_{nC} = \cc_{np} \bveps_{pG} + \cc_{nn} \bveps_{nG}
\label{eqn:stressdef}
\end{eqnarray}
with
\begin{eqnarray}
\cc_{pp} = \left[ \begin{array}{ccc} C_{11} & C_{12} & C_{16} \\ C_{12} & C_{22} & C_{26} \\ C_{16} & C_{26} & C_{66} \end{array} \right] \hspace{0.5cm} \cc_{pn} = \left[ \begin{array}{ccc} 0 & 0 & C_{13} \\ 0 & 0 & C_{23} \\ 0 & 0 & C_{36} \end{array} \right] \nonumber \\
\cc_{np} = \left[ \begin{array}{ccc} 0 & 0 & 0 \\ 0 & 0 & 0 \\ C_{13} & C_{23} & C_{36} \end{array} \right] \hspace{0.5cm} \cc_{nn} = \left[ \begin{array}{ccc} C_{55} & C_{45} & 0 \\ C_{45} & C_{44} & 0 \\ 0 & 0 & C_{33} \end{array} \right]
\end{eqnarray}
where the subscript $C$ indicate the constitutive equations. The \textit{Principle of Virtual Displacements} (PVD) in case of multilayered plate subjected to mechanical loads is written as:
\begin{equation}
\sum\limits_{k=1}^{N_k} \int\limits_{\Omega_k}\int\limits_{A_k} \left\{  (\delta \bveps_{pG}^k)^{\rm T} \bvsig_{pC}^k + (\delta \bveps_{nG}^k)^{\rm T} \bvsig_{nC}^k \right\}~\mathrm{d}\Omega_k~\mathrm{d}z = \sum\limits_{k=1}^{N_k} \int\limits_{\Omega_k} \int\limits_{A_k} \rho^k \delta \uu_s^{k^{\rm T}} \ddot{\uu}^k~\mathrm{d}\Omega_k~\mathrm{d}z + \sum\limits_{k=1}^{N_k} \delta \mathbf{L}_e^k
\end{equation}
where $\rho^k$ is the mass density of the $k^{\rm th}$ layer, $\Omega_k$, $A_k$ are the integration domain in the $(x,y)$ and the $z$ direction, respectively. Upon substituting the geometric relations (\Eref{eqn:strainDisp}), the constitutive relations (\Eref{eqn:stressdef}) and the unified formulation into the PVD statement, we have:
\begin{equation}
\begin{split}
\int\limits_{\Omega_k}\int\limits_{A_k} \left\{ \left(\dd_p^k F_s \delta \uu_s^k\right)^{\rm T} \left\{ \cc_{pp}^k \dd_p^k F_{\tau} \uu_{\tau}^k + \cc_{pn}^k (\dd_{n\Omega}^k + \dd_{nz}^k) F_{\tau}\uu_\tau^k \right\} + \right. \\
\left. \left[ (\dd_{n\Omega}^k + \dd_{nz}^k) f_x \delta\uu_s^k)^{\rm T} (\cc_{np}^k \dd_p^k F_\tau \uu_\tau^k  + \cc_{nn}^k (\dd_{n\Omega}^k + \dd_{nz}^k) F_\tau \uu_\tau^k) \right] \right\}~\mathrm{d}\Omega_k ~\mathrm{d}z = \\ \sum\limits_{k=1}^{N_k} \int\limits_{\Omega_k} \int\limits_{A_k} \rho^k \delta \uu_s^{k^{\rm T}} \ddot{\uu}^k~\mathrm{d}\Omega_k~\mathrm{d}z + \sum\limits_{k=1}^{N_k} \delta \mathbf{L}_e^k
\end{split}
\end{equation}
After integration by parts, the governing equations for the plate are obtained:
\begin{equation}
\kk_{uu}^{k\tau s} \uu_{\tau}^k = \mathbf{P}_{u \tau}^k
\end{equation}
and in the case of free vibrations, we have:
\begin{equation}
\kk_{uu}^{k\tau s} \uu_{\tau}^k = \mathbf{M}^{k \tau s} \ddot{\uu}_\tau^k
\end{equation}
where the fundamental nucleus $\kk_{uu}^{k \tau s}$ is:
\begin{equation}
\kk_{uu}^{k \tau s} = \left[ (-\dd_p^k)^{\rm T} ( \cc_{pp}^k \dd_p^k + \cc_{pn}^k (\dd_{n\Omega}^k + \dd_{nz}) + (-\dd_{n\Omega}^k + \dd_{nz}^k)^{\rm T} (\cc_{np}^k \dd_p^k + \cc_{nn}^k (\dd_{n\Omega}^k + \dd_{nz}^k)) \right] F_\tau F_s
\label{eqn:stifffundanuclei}
\end{equation}
and $\mathbf{M}^{k \tau s}$ is the fundamental nucleus for the inertial term given by:
\begin{equation}
M_{ij}^{k \tau s} = \left\{ \begin{array}{cc} \rho^k F_\tau F_s & \textup{if} \hspace{1cm} i = j \\ 0 & \textup{if} \hspace{1cm} i \neq j \end{array} \right.
\label{eqn:massfundanuclei}
\end{equation}
where $\mathbf{P}_{u \tau}^k$ are variationally consistent loads with applied pressure. For more detailed derivation and for the explicit form of the fundamental nuclei, interested readers are referred to~\cite{carrerademasi2002,carrerademasi2002a}.

\section{Non-uniform rational B-splines}\label{isog}
In this study, the finite element approximation uses NURBS basis function. We give here only a brief introduction to NURBS. More details on their use in FEM are given in~\cite{cottrellhughes2009}. The key ingredients in the construction of NURBS basis functions are: the knot vector (a non decreasing sequence of parameter values, $\xi_i \le \xi_{i+1}, i = 0,1,\cdots,m-1$), the control points, $P_i$, the degree of the curve $p$ and the weight associated to a control point, $w$. The i$^{th}$ B-spline basis function of degree $p$, denoted by $N_{i,p}$ is defined as: 

\begin{eqnarray}
N_{i,0}(\xi) = \left\{ \begin{array}{cc} 1 & \textup{if} \hspace{0.2cm} \xi_i \le \xi \le \xi_{i+1} \\
0 & \textup{else} \end{array} \right. \nonumber \\
N_{i,p}(\xi) = \frac{ \xi- \xi_i}{\xi_{i+p} - \xi_i} N_{i,p-1}(\xi) + \frac{\xi_{i+p+1} - \xi}{\xi_{i+p+1}-\xi_{i+1}}N_{i+1,p-1}(\xi)
\end{eqnarray}
A $p^{th}$ degree NURBS curve is defined as follows:

\begin{equation}
\mathbf{C}(\xi) = \frac{\sum\limits_{i=0}^m N_{i,p}(\xi)w_i \mathbf{P}_i} {\sum\limits_{i=0}^m N_{i,p}(\xi)w_i}
\label{eqn:nurbsfunc1}
\end{equation}

\begin{figure}[htpb]
\centering
\includegraphics[scale=0.6]{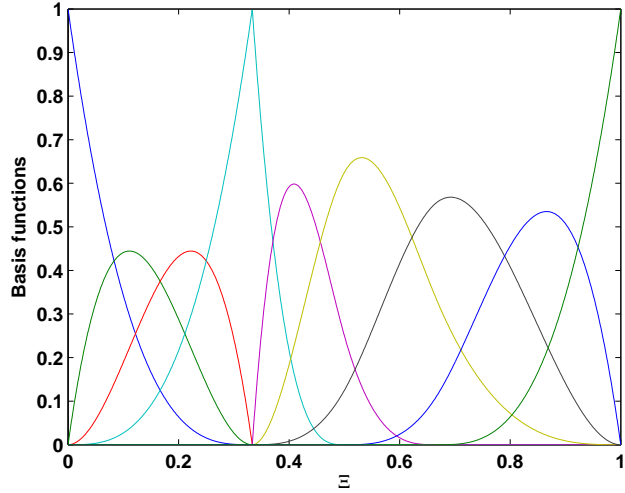}
\caption{non-uniform rational B-splines, order of the curve = 3}
\label{fig:nurbsplot}
\end{figure}
where $\mathbf{P}_i$ are the control points and $w_i$ are the associated weights. \fref{fig:nurbsplot} shows the third order non-uniform rational B-splines for a knot vector, $\Xi = \{0,~ 0,~ 0,~ 0,~ 1/3,~ 1/3,~ 1/3,~ 1/2,~ 2/3,~ 1,~ 1,~ 1,~ 1\}$. NURBS basis functions has the following properties: (i) non-negativity, (ii) partition of unity, $\sum\limits_i N_{i,p} = 1$; (iii) interpolatory at the end points. As the same function is also used to represent the geometry, the exact representation of the geometry is preserved. It should be noted that the continuity of the NURBS functions can be tailored to the needs of the problem. The B-spline surfaces are defined by the tensor product of basis functions in two parametric dimensions $\xi$ and $\eta$ with two knot vectors, one in each dimension as:
\begin{equation}
\mathbf{C}(\xi,\eta) = \sum_{i=1}^n\sum_{j=1}^m N_{i,p}(\xi)M_{j,q}(\eta) \mathbf{P}_{i,j}
\end{equation}
where $\mathbf{P}_{i,j}$ is the bidirectional control net and $N_{i,p}$ and $M_{j,q}$ are the B-spline basis functions defined on the knot vectors over an $m\times n$ net of control points $\mathbf{P}_{i,j}$. The NURBS surface is then defined by:
\begin{equation}
\mathbf{C}(\xi,\eta) = \frac{\sum_{i=1}^n\sum_{j=1}^m N_{i,p}(\xi)M_{j,q}(\eta) \mathbf{P}_{i,j} w_iw_j}{w(\xi,\eta)}
\label{eqn:nurbsfunc}
\end{equation}
where $w(\xi,\eta)$ is the weighting function. The displacement field within the control mesh is approximated by:
\begin{equation}
\uu_\tau(x,y) = \mathbf{C}(\xi,\eta) \qq_\tau(x,y),
\end{equation}
where $\qq_\tau(x,y)$ are the nodal variables and $\mathbf{C}(\xi,\eta)$ are the basis functions given by \Eref{eqn:nurbsfunc}.

\subsection{Shear locking} Similar to the finite element based on Lagrange basis functions, locking appears when lower order NURBS basis functions are employed~\cite{veigabuffa2012,valizadehnatarajan2013}, for example with quadratic, cubic and quartic elements~\footnote{Linear NURBS basis functions are same as the linear Lagrange basis functions and are not discussed here. Approaches employed for Lagrange basis functions can readily be applied to NURBS basis functions with order 1}. One approach to alleviate the shear locking is to employ interpolation functions of order 5 or higher~\cite{veigabuffa2012}, but this inevitably increases the computational cost. A stabilization technique for several lower-order NURBS elements for plates was reported in \cite{ChienHungNhon2012}. In this paper, we adopt a stabilization technique proposed in~\cite{kikuchiishii1999} and later used in~\cite{valizadehnatarajan2013} to study the response of Reissner-Mindlin plates. In this approach, the material matrix related to the shear terms are multiplied by the following factor:
\begin{equation}
\textup{shearFactor} = \frac{h^2}{h^2 + \alpha^2 \ell^2}
\end{equation}
where $\ell$ is the longest length of the edges of the NURBS element and $\alpha$ is a positive constant given in the interval $0.05 \le \alpha \le 0.15$. It is found from numerical experiments of NURBS-based isogeometric plate elements that $\alpha$ can be fixed at 0.1, which provide reasonably accurate solutions.

\section{Numerical Results}
\label{numres}
In this section, we present the static response and the natural frequencies of laminated composite plates using the combined IGA and CUF framework. In this study we use a hybrid displacement assumption, where the in-plane displacements $u$ and $v$ are expressed as sinusoidal expansion in the thickness direction, and the transverse displacement, $w$ is quadratic in the thickness direction. We refer to this theory as SINUS-W2. The displacements are expressed as:
\begin{align}
u(x,y,z,t) &= u_o(x,y,t) + zu_1(x,y,t) + \sin \left( \frac{\pi z}{h} \right) u_2(x,y,t) \nonumber \\
v(x,y,z,t) &= v_o(x,y,t) + zv_1(x,y,t) + \sin \left( \frac{\pi z}{h}\right) v_2(x,y,t) \nonumber \\
w(x,y,z,t) &= w_o(x,y,t) + zw_1(x,y,t) + z^2 w_2(x,y,t)
\end{align}
where $u_o, v_o$ and $w_o$ are translations of a point at the middle-surface of the plate, $w_2$ is higher order translation, and $u_1, v_1, u_3$ and $v_3$ denote rotations~\cite{touratier1991} and considers a quadratic variation of the transverse displacement $w$ allowing for the through-the-thickness deformations. The effect of the plate aspect ratio, the ply angle and the ratio of Young's modulus $E_1/E_2$ on the static bending and free vibration is numerically studied.

\subsection{Static bending} The static analysis is conducted for cross-ply laminated plates with three and four layers under the following sinusoidal load:
\begin{equation}
p_z(x,y) = P_o \sin \left( \frac{\pi x}{a} \right) \sin \left( \frac{\pi y}{a} \right)
\label{eqn:mechload}
\end{equation}
where $P_o$ is the amplitude of the mechanical load. The origin of the coordinate system is located at the lower-left corner on the midplane. The physical quantities are non-dimensionalized by the following relations, unless otherwise mentioned:
\begin{eqnarray}
\overline{w} = w(a/2,a/2,0)\frac{100 h^3E_2}{Pa^4}; ~~ \overline{\sigma}_{xx} = \sigma_{xx}(a/2,a/2,h/2) \frac{h^2}{Pa^2}; \nonumber \\
\overline{\sigma}_{yy} = \sigma_{yy}(a/2,a/2,h/4)\frac{h^2}{Pa^2}; ~~ \overline{\tau}_{xz} = \tau_{xz}(0,a/2,0) \frac{h}{Pa};
\end{eqnarray}

\paragraph{Validation} Before proceeding with a detailed numerical study on the effect of various parameters on the global response of cross-ply laminated composites, the results from the proposed formulation are compared against available results pertaining to static bending of laminated plates. In this study, we consider three orders of NURBS basis functions, viz., quadratic, cubic and quartic. It is noted that, in this study, we do not consider first order NURBS basis functions. This is because, the first order NURBS basis functions are similar to the conventional bilinear shape functions. The performance of which is discussed in detail in~\cite{carreracinefra2010,natarajanferreira2013}. In this study, the results from the present formulation are denoted by Quadratic, Cubic and Quartic, which corresponds to the order of shape functions employed, which is referred to as $p-$refinement. Three different mesh discretizations, viz., 5$\times$5, 7$\times$7 and 9$\times$9 are considered, which is called as $h-$refinement. Table \ref{tab:meshconvedeflec} shows the convergence of the central deflection and stresses of a simply supported cross-ply laminated square plate. It is seen that with both $h-$ and $p-$refinement, the results from the present formulation converge. It is seen that highly accurate results are obtained from the present formulation even with a coarse mesh. A comparison with other approaches and an elasticity solution is given in Table \ref{tab:fourlayerstaticdeflec}.

\begin{table}[htbp]
\centering
\renewcommand{\arraystretch}{1.2}
\caption{Convergence of the central deflection $\overline{w} = w(a/2,a/2,0)\frac{100 E_2 h^3}{P a^4}$ of a simply supported cross-ply laminated square plate $[0^\circ/90^\circ/90^\circ/0^\circ]$ with $E_1=$ 25$E_2$, $G_{12}=G_{13}=$ 0.5$E_2$, $G_{23}=$ 0.2$E_2$, $\nu_{12}=$0.25. }
\begin{tabular}{llcrrr}
\hline
 &  Method && \multicolumn{3}{c}{Meshes} \\
\cline{3-6}
& && 5$\times$5 & 7$\times$ 7 & 9$\times$9  \\
\hline
\multirow{4}{*}{$\overline{w}$} & Quadratic &&  1.9207 & 1.9100 & 1.9058\\
& Cubic &&  1.9076 & 1.9038 & 1.9021 \\
& Quartic && 1.9045 & 1.9020 & 1.9010 \\
& HSDT~\cite{reddy1984} && 1.8937  \\
& Elasticity~\cite{pagano1970}  && 1.9540 \\
\cline{2-6}
\multirow{4}{*}{$\overline{\sigma}_{xx}$} & Quadratic && 0.6966 & 0.7009 & 0.7029 \\
& Cubic &&   0.7074 & 0.7063	 & 0.7061 \\
& Quartic && 0.7062 & 0.7060 & 0.7058 \\
& HSDT~\cite{reddy1984} && 0.6651 \\
& Elasticity~\cite{pagano1970}  && 0.7200 \\
\cline{2-6}
\multirow{4}{*}{$\overline{\sigma}_{yy}$} & Quadratic && 0.6179 & 0.6221 & 0.6239\\
& Cubic &&   0.6277 & 0.6270	 & 0.6268 \\
& Quartic && 0.6268 & 0.6267 & 0.6266 \\
& HSDT~\cite{reddy1984} && 0.6322 \\
& Elasticity~\cite{pagano1970}  && 0.6660 \\
\cline{2-6}
\multirow{4}{*}{$\overline{\tau}_{xz}$} & Quadratic && 0.2293 & 0.2246 & 0.2227 \\
& Cubic && 0.2210 & 0.2205 & 0.2202 \\
& Quartic && 0.2205 & 0.2202 & 0.2201 \\
& HSDT~\cite{reddy1984} && 0.2064 \\
& Elasticity~\cite{pagano1970}  && 0.2700 \\
\cline{2-6}
\hline
\end{tabular}
\label{tab:meshconvedeflec}
\end{table}

\subsubsection{Four layer (0$^\circ$/90$^\circ$)$_{\rm s}$ square cross-ply laminated plate under sinusoidal load} A square simply supported laminate of side $a$ and thickness $h$, composed of four equally thick layers oriented at (0$^\circ$/90$^\circ$)$_{\rm s}$ is considered. The plate is subjected to a sinusoidal vertical pressure given by \Eref{eqn:mechload}. The material properties are as follows: $E_1= \textup{25}E_2;~G_{12}=G_{13} = \textup{0.5}E_2;~G_{23}=\textup{0.2}E_2;~\nu_{12}=\textup{0.25}$. For this example, a three-dimensional exact solution by Pagano~\cite{pagano1970} is available. The central deflection and the corresponding stresses for the SINUS-W2 theory with isogeometric approach are presented in Table \ref{tab:fourlayerstaticdeflec}. We compare the results with higher order plate theories~\cite{reddy1984,ferreiracarrera2012}, first order theory~\cite{reddychao1981}, an exact solution~\cite{pagano1970} and also with the strain smoothing approach with SINUS-W2. The effect of plate the thickness is also shown in Table \ref{tab:fourlayerstaticdeflec}. It is clear that the first order shear deformation theories (FSDT) cannot be used for thick laminates. It can be seen that the results from the present formulation are in very good agreement with those in the literature and very precise transverse displacements and stresses are obtained. 

\begin{table}[htbp]
\centering
\renewcommand{\arraystretch}{1.5}
\caption{The normalized central deflection $\overline{w} = w(a/2,a/2,0)\frac{100 E_2 h^3}{P a^4}$, stresses, $\overline{\sigma}_{xx} = \sigma_{xx}(a/2,a/2,h/2)\frac{h^2}{Pa^2}, \overline{\sigma}_{yy} = \sigma_{yy}(a/2,a/2,h/4) \frac{h^2}{Pa^2}$ and $\overline{\tau}_{xz} = \tau_{xz}(0,/a2,0) \frac{h}{Pa}$  of a simply supported cross-ply laminated square plate $[0^\circ/90^\circ/90^\circ/0^\circ]$, with $E_1=$ 25$E_2$, $G_{12}=G_{13}=$ 0.5$E_2$, $G_{23}=$ 0.2$E_2$, $\nu_{12}=$0.25. }
\begin{tabular}{llcccc}
\hline
$a/h$ & Method & $w$ & $\sigma_{xx}$ & $\sigma_{yy}$ & $\tau_{xz}$\\
\hline
\multirow{5}{*}{10} & HSDT~\cite{reddy1984} & 0.7147 & 0.5456 & 0.3888 & 0.2640 \\
& FSDT~\cite{reddychao1981} & 0.6628 & 0.4989 & 0.3615 & 0.1667 \\
& Elasticity~\cite{pagano1970} & 0.7430 & 0.5590 & 0.4030 & 0.3010\\
& RBF~\cite{ferreiracarrera2012} & 0.7325 & 0.5627 & 0.3908 & 0.3321 \\
& CS-FEM Q4 (4 subcells)~\cite{natarajanferreira2013} & 0.7195 & 0.5597 & 0.3905 & 0.2952 \\
& Present (Quadratic 9$\times$ 9) & 0.7250 & 0.5571 & 0.3908 & 0.2985\\
& Present (Cubic 9$\times$9) & 0.7203 & 0.5596 & 0.3913 & 0.2983 \\
& Present (Quartic 9$\times$9) & 0.7187 & 0.5594 & 0.3907 & 0.2967 \\
\cline{2-6}
\multirow{5}{*}{100} & HSDT~\cite{reddy1984} & 0.4343 & 0.5387 & 0.2708 & 0.2897\\
& FSDT~\cite{reddychao1981} & 04337 & 0.5382 & 0.2705 & 0.1780 \\
& Elasticity~\cite{pagano1970} & 0.4347 & 0.5390 & 0.2710 & 0.3390\\
& RBF~\cite{ferreiracarrera2012} & 0.4307 & 0.5431 & 0.2730 & 0.3768\\
& CS-FEM Q4 (4 subcells)~\cite{natarajanferreira2013} & 0.4304 & 0.5368 & - & 0.3285 \\
& Present (Quadratic 9$\times$ 9) & 0.4383 & 0.5334 & - & 0.4069 \\
& Present (Cubic 9$\times$9) & 0.4336 & 0.5368 & - &	0.3271 \\
& Present (Quartic 9$\times$9) & 0.4317 & 0.5366 & - & 0.3275 \\
\hline
\end{tabular}
\label{tab:fourlayerstaticdeflec}
\end{table}

\subsubsection{Three layer (0$^\circ$/90$^\circ$/0$^\circ$) square cross ply laminated plate under sinusoidal load} In this case, a square laminate of side $a$ and thickness $h$, composed of three equally thick layers oriented at (0$^\circ$/90$^\circ$/0$^\circ$) is considered. It is simply supported on all edges and subjected to a sinusoidal vertical pressure of the form given by \Eref{eqn:mechload}. The material properties for this example are: $E_1=$132.38 GPa, $E_2=E_3=$10.756 GPa, $G_{12}=$3.606 GPa, $G_{13}=G_{23}=$ 5.6537 GPa, $\nu_{12}=\nu_{13}=$ 0.24, $\nu_{23}=$ 0.49. In Table \ref{tab:threelayerstaticdeflec}, we present results for the SINUS-W2 theory with isogeometric approach with quadratic, cubic and quartic NURBS basis functions with a 9$\times$9 NURBS patch. The results from the present formulation are compared with the analytical solution~\cite{carrera1998,carrera2001} and the MITC4 formulation with and without strain smoothing~\cite{carreracinefra2010,natarajanferreira2013}. It can be seen that the numerical results from the present formulation are found to be in good agreement with the existing solutions. Moreover, it is noted that with the isogeometric approach, the geometry of the domain can be exactly represented. Although only simple geometry is considered, the proposed formulation easily be extended to complex geometry. The main features of the present formulation are: (1) theories from ESL to higher order layer descriptions can be implemented within a single code (since it is based on CUF); (2) the isogeometric approach provides flexibility to construct higher order smooth functions and provides accurate solutions even for a coarse NURBS mesh and (3) the present formulation is insensitive to shear locking.
\begin{table}[htpb]
\centering
\renewcommand{\arraystretch}{1.5}
\caption{Transverse displacement $\overline{w} = w(a/2,a/2,h/2)$ at the center of a multilayered plate $[0^\circ/90^\circ/0^\circ]$ with $E_1=$ 132.38 GPa, $E_2=E_3=$ 10.756 GPa, $G_{12}=$ 3.606 GPa, $G_{13}=G_{23}=$ 5.6537 GPa, $\nu_{12}=\nu_{13}=$ 0.24, $\nu_{23}=$ 0.49. }
\begin{tabular}{lccccc}
\hline
$\overline{w}$ & \multicolumn{5}{c}{$a/h$} \\
\cline{2-6}
& 10 & 50 & 100 & 500 & 1000 \\
\hline
Analytical (ESL-2)~\cite{carrera1998,carrera2001} & 0.9249 & 0.7767 & 0.7720 & 0.7705 & 0.7704 \\
MITC4~\cite{carreracinefra2010} & 0.9195 & 0.7713 & 0.7666 & 0.7650 & 0.7650 \\
CS-FEM Q4 (4 subcells)~\cite{natarajanferreira2013} & 0.9235 & 0.7703 & 0.7655 & 0.7639 & 0.7639\\
Present (Quadratic 9$\times$9) & 0.9252 & 0.7713 & 0.7650 & 0.7624 & 0.7624\\
Present (Cubic 9$\times$9) & 0.9226 & 0.7704 & 0.7656 & 0.7640 & 0.7639 \\
Present (Quartic 9$\times$9) & 0.9217 & 0.7695 & 0.7646 & 0.7631 & 0.7630\\
\hline
\end{tabular}
\label{tab:threelayerstaticdeflec}
\end{table}

\subsection{Free vibration - cross-ply laminated plates}
In this example, all layers of the laminate are assumed to be of the same thickness, density and made up of the same linear elastic material. The following material parameters are considered for each layer
\begin{eqnarray*}
\frac{E_1}{E_2} = \textup{10,20,30,~or~ 40}; ~~ G_{12}=G_{13} = \textup{0.6}E_2; \\
G_3 = \textup{0.5}E_2; \nu_{12} = \textup{0.25}.
\end{eqnarray*}
The subscripts 1 and 2 denote the directions normal and the transverse to the fiber direction in a lamina, which may be oriented at an angle to the plate axes. The ply angle of each layer is measure from the global $x-$axis to the fiber direction. The example considered is a simply supported square plate of the cross-ply lamination  (0$^\circ$/90$^\circ$)$_{\rm s}$. The thickness and the length of the plate are denoted by $h$ and $a$, respectively. The thickness-to-span ratio $h/a=$ 0.2 is employed in the computations. In this study, we present the non dimensionalized free flexural frequencies as, unless specified otherwise:
\begin{equation*}
\Omega = \omega \frac{a^2}{h}\sqrt{ \frac{\rho}{E_2} }
\end{equation*}

\begin{table}[htbp]
\centering
\renewcommand{\arraystretch}{1.5}
\caption{Convergence of the normalized fundamental frequency $\Omega = \omega a^2/h\sqrt{\rho/E_2}$ of a simply supported cross-ply laminated square plate  (0$^\circ$/90$^\circ$)$_{\rm s}$  with $h/a=0.2$, $\frac{E_1}{E_2}=$ 40, $G_{12}=G_{13}=$ 0.6$E_2$, $G_{23}=$ 0.5$E_2$, ~$\nu_{12}=$ 0.25.. }
\begin{tabular}{lccc}
\hline
Method & \multicolumn{3}{c}{Meshes} \\
\cline{2-4}
& 5$\times$5 & 7$\times$ 7 & 9$\times$9  \\
\hline
Quadratic & 10.6926 & 10.7295 & 10.7454 \\
Cubic & 10.7340 & 10.7517 & 10.7590 \\
Quartic & 10.7498 & 10.7598 & 10.7640\\
\hline
\end{tabular}
\label{tab:meshconveFreq}
\end{table}

\begin{table}[htbp]
\centering
\renewcommand{\arraystretch}{1.2}
\caption{The normalized fundamental frequency $\Omega = \omega a^2/h\sqrt{\rho/E_2}$ of a simply supported cross-ply laminated square plate  (0$^\circ$/90$^\circ$)$_{\rm s}$  with $h/a=0.2$, $\frac{E_1}{E_2}=$ 10, 20, 30 or 40, $G_{12}=G_{13}=$ 0.6$E_2$, $G_{23}=$ 0.5$E_2$, ~$\nu_{12}=$ 0.25.}
\begin{tabular}{lccccc}
\hline
Method & \multicolumn{4}{c}{$E_1/E_2$} \\
\cline{2-5}
 &    10 & 20 & 30 & 40 \\
\hline
Liew~\cite{liewhuang2003} &  8.2924 & 9.5613 & 10.3200 & 10.8490 \\
Reddy, Khdeir~\cite{khdeirlibrescu1988}   & 8.2982 & 9.5671 & 10.3260 & 10.8540 \\
HSDT~\cite{ferreiraroque2011} $(\nu_{23}=0.18)$ & 8.2999 & 9.5411 & 10.2687 & 10.7652 \\
CS-FEM Q4 (4 subcells)~\cite{natarajanferreira2013}  &  8.3642 & 9.5793 & 10.2973 & 10.7887\\
Present (Quadratic 9$\times$9) & 8.3358 & 9.5437 & 10.2572 & 10.7454 \\
Present (Cubic 9$\times$9) & 8.3417 & 9.5532 & 10.2691 & 10.7590 \\
Present (Quartic 9$\times$9) & 8.3439 & 9.5566 & 10.2734 & 10.7640 \\
\hline
\end{tabular}
\label{tab:freevib1}
\end{table}
Table \ref{tab:meshconveFreq} shows the convergence of the normalized fundamental frequency of a simply supported cross-ply laminated square plate based on the current isogeometric approach. The performance of various basis functions with NURBS mesh refinement is studied. It is seen that with $h-$ refinement, the solutions converge and with $p-$ refinement, the accuracy increases for the same mesh size, as expected. Table \ref{tab:freevib1} lists the fundamental frequency for a simply supported cross-ply laminated square plate with $h/a=$ 0.2 and for different Young's modulus ratios, $E_1/E_2$. It can be seen that the results from the present formulation are in very close agreement with the values of \cite{khdeirlibrescu1988} based on higher order theory, the meshfree results of Liew \textit{et al.,}~\cite{liewhuang2003} and Ferreira \textit{et al.,} based on FSDT and higher order theories with radial basis functions~\cite{ferreiraroque2011}. The effect of plate thickness on the fundamental frequency is shown in Table \ref{tab:freevib2}. It can be seen that the results agree with the results available in the literature. The present formulation is insensitive to shear locking.

\begin{table}[htbp]
\centering
\renewcommand{\arraystretch}{1}
\caption{Variation of fundamental frequencies, $\Omega = \omega a^2/h\sqrt{\rho/E_2}$ with $a/h$ for a simply supported square laminated plate $[0^\circ/90^\circ/90^\circ/0^\circ]$, $\Omega = \omega a^2/h\sqrt{\rho/E_2}$,  with $E_1/E_2=$ 40, $G_{12}=G_{13}=$0.6$E_2$, $G_{23}=$0.5$E_2$, $\nu_{12}=\nu_{13}=\nu_{23}=$ 0.25.}
\begin{tabular}{lcccccc}
\hline
Method & \multicolumn{6}{c}{$a/h$}\\
\cline{2-7}
 &  2 & 4 & 10 & 20 & 50 & 100 \\
\hline
FSDT~\cite{whitneypagano1970} & 5.4998 & 9.3949 & 15.1426 & 17.6596 & 18.6742 & 18.8362 \\
Model-1 (12dofs)~\cite{kantswaminathan2001} & 5.4033 & 9.2870 & 15.1048 & 17.6470 & 18.6720 & 18.8357 \\
Model-2 (9dofs)~\cite{kantswaminathan2001} & 5.3929 & 9.2710 & 15.0949 & 17.6434 & 18.6713 & 18.8355 \\
HSDT~\cite{reddy1984} & 5.5065 & 9.3235 & 15.1073 & 17.6457 & 18.6718 & 18.8356 \\
HSDT~\cite{senthilnathanlim1987}& 6.0017 & 10.2032 & 15.9405 & 17.9938 & 18.7381 & 18.8526 \\
CS-FEM Q4 (4 subcells)~\cite{natarajanferreira2013} & 5.4026 & 9.2998 & 15.1766 & 17.7540 & 18.7947 & 18.9611 \\
Present (Quadratic 9$\times$9) &5.3931 & 9.2701 & 15.0660 & 17.5781 & 18.5913 & 18.7579 \\
Present (Cubic 9$\times$9) &5.3945 & 9.2785 & 15.1086 & 17.649 & 18.6711 & 18.8343 \\
Present (Quartic 9$\times$9) & 5.3951 & 9.2815 & 15.1239 & 17.6749 & 18.7024 & 18.8665 \\
\hline
\end{tabular}
\label{tab:freevib2}
\end{table}

\subsection{Circular plates}
In this example, consider a circular four layer $[\theta/-\theta/-\theta/\theta]$ laminated plate with fully clamped boundary conditions. The influence of the fiber orientations on the free vibration of clamped circular laminated plate is studied. The following material properties are used:
\begin{eqnarray*}
\frac{E_1}{E_2} = \textup{40}; ~~ G_{12}=G_{13} = \textup{0.6}E_2; \\
G_3 = \textup{0.5}E_2; \nu_{12} = \textup{0.25}.
\end{eqnarray*}
The subscripts 1 and 2 denote the directions normal and the transverse to the fiber direction in a lamina. The circular plate has a radius-to-thickness of 5 ($R/h=$5). For this problem, a NURBS quadratic basis function is enough to model exactly the circular geometry. Any further refinement, if done, will only improve the accuracy of the solution. The following knot vectors for the coarsest mesh with one element are defined as follows: $\Xi=$[0,0,0,1,1,1]; and $\mathcal{H}=$[0,0,0,1,1,1]. The data for the circular plate is given in Table \ref{table:circplatedata}. In this study, 13$\times$13 NURBS cubic elements are used. The first three fundamental frequencies for a clamped circular laminated plate are given in Table \ref{table:circplatefreq}. The fiber-orientation of each layer is considered to be the same and the influence of the fiber orientation on the first three fundamental frequencies are given in Table \ref{table:circplatefreq}. The numerical results from the present approach are compared with the moving least square differential quadrature method (MLSDQ) based on FSDT~\cite{liewhuang2003} and IGA with inverse trigonometric shear deformation theory~\cite{thaiferreira2014}. It can be seen that the results from the present formulation agree well with the results in the literature.

\begin{table}[htbp]
\centering
\renewcommand{\arraystretch}{1}
\caption{Control points and weights for a circular plate with radius $R=$ 0.5.}
\begin{tabular}{lrrrrrrrrr}
\hline
i & 1 & 2 & 3 & 4 & 5 & 6 & 7 & 8 & 9 \\
\hline
$x_i$ & -$\frac{\sqrt{2}}{4}$ & -$\frac{\sqrt{2}}{2}$ & $\frac{\sqrt{2}}{4}$ & 0 & 0 & 0 & $\frac{\sqrt{2}}{4}$ & $\frac{\sqrt{2}}{2}$ & $\frac{\sqrt{2}}{4}$ \\
$y_i$ & $\frac{\sqrt{2}}{4}$ & 0 & -$\frac{\sqrt{2}}{4}$ & $\frac{\sqrt{2}}{2}$ & 0 & -$\frac{\sqrt{2}}{2}$ & $\frac{\sqrt{2}}{4}$ & 0 &-$\frac{\sqrt{2}}{4}$ \\ 
$w_i$ & 1 & $\frac{\sqrt{2}}{2}$ & 1 & $\frac{\sqrt{2}}{2}$ & 1 & $\frac{\sqrt{2}}{2}$ & 1 & $\frac{\sqrt{2}}{2}$ & 1 \\
\hline
\end{tabular}
\label{table:circplatedata}
\end{table}

\begin{table}[htbp]
\centering
\renewcommand{\arraystretch}{1}
\caption{Influence of fiber orientations on the fundamental frequencies, $\Omega=\omega a^2/h \sqrt{\rho/E_2}$ for clamped circular laminated plates.}
\begin{tabular}{llcrrr}
\hline
$\theta$ & Method & & \multicolumn{3}{c}{$\Omega$}\\
\cline{4-6}
& && 1 & 2 & 3 \\
\hline
\multirow{3}{*}{0} & MLSDQ-FSDT~\cite{liewhuang2003} && 22.2110 & 29.651 & 41.1010 \\
& IGA~\cite{thaiferreira2014} & & 23.5781 & 30.7459 & 42.0042 \\
& Present && 22.6663 & 30.3485 & 41.7294 \\
\cline{2-6}
\multirow{3}{*}{$\pi/12$} & MLSDQ-FSDT~\cite{liewhuang2003} && 22.7740 & 31.4550 & 43.350 \\
& IGA~\cite{thaiferreira2014} && 23.6090 & 31.7743 & 43.9569 \\
& Present && 23.0024 & 31.5752 & 43.7671 \\
\cline{2-6}
\multirow{3}{*}{$\pi/6$} & MLSDQ-FSDT~\cite{liewhuang2003} && 24.0710 & 36.1530 & 43.9680 \\
& IGA~\cite{thaiferreira2014} && 24.2081 & 35.6047 & 46.5406 \\
& Present && 23.9749 & 35.2577 & 44.2964 \\
\cline{2-6}
\multirow{3}{*}{$\pi/4$} & MLSDQ-FSDT~\cite{liewhuang2003} && 24.7520 & 39.1810 & 43.6070 \\
& IGA~\cite{thaiferreira2014} && 24.6607 & 37.8980 & 46.2560 \\
& Present && 24.5253 & 37.4311 & 44.0796 \\
\hline
\end{tabular}
\label{table:circplatefreq}
\end{table}


\section{Conclusions}
In this article, the isogeometric approach was combined with the unified formulation to study the static bending and the free vibration of laminated composites. The present approach allows us to achieve smooth approximation of the unknown fields with arbitrary continuity. When employing lower order elements, the method suffers from shear locking syndrome, which is alleviated by multiplying the shear term with a correction factor. The results from the present formulations are in very good agreement with the solutions available in the literature. It is believed that the present formulation is definitely a effective computational formulation for practical problems. On one hand, the unified formulation allows the user to test different theories within a single framework, whilst, the isogeometric approach not only provides flexibility in constructing higher smooth basis functions, but the geometry is accurately described.

\section*{Acknowledgements}
S Natarajan would like to acknowledge the financial support of the School of Civil and Environmental Engineering, The University of New South Wales for his research fellowship for the period September 2012 onwards.

\bibliographystyle{elsarticle-num}
\bibliography{laminatedCompo}

\end{document}